\input amstex
\input amsppt.sty
\magnification\magstep1

\def\ni\noindent
\def\sbs{\subset}

\def\diam{\operatorname{diam}}

\def\R{\text{\bf R}}

\def\Z{\text{\bf Z}}
\def\H{\text{\bf H}}

\def\N{\text{\bf N}}

\def\G{\Gamma}

\hoffset= 0.0in
\voffset= 0.0in
\hsize=32pc
\vsize=38pc
\baselineskip=24pt
\NoBlackBoxes
\topmatter
\author
Alexander Dranishnikov and Viktor Schroeder
\endauthor

\title
Aperiodic colorings and tilings of Coxeter groups
\endtitle
\abstract
We construct a limit aperiodic coloring of hyperbolic groups.
Also we construct
limit strongly aperiodic strictly balanced tilings
of the Davis complex for all Coxeter groups.
\endabstract

\thanks The authors were partially supported by NSF grants DMS-0604494
and SNF grant 200020-103594
\endthanks

\address University of Florida, Department of Mathematics, P.O.~Box~118105,
358 Little Hall, Gainesville, FL 32611-8105, USA;
University of Zurich, Department of Mathematics,
Winterthurerstrasse 190, CH-8057 Zurich, Switzerland
\endaddress


\email  dranish\@math.ufl.edu;
vschroed\@math.unizh.ch
\endemail

\endtopmatter

\document
\head \S0 Introduction \endhead

In [BDS] we constructed a quasi-isometric embedding of hyperbolic
groups into a finite product of binary trees. First we implemented
such construction for hyperbolic Coxeter groups [DS]. As a
byproduct we obtained aperiodic tilings with finitely many tiles
of the Davis complex for these groups. Our tilings are limit
strongly aperiodic and the set of tiles can be taken to be
aperiodic. As result we obtain limit strongly aperiodic
tilings of the hyperbolic spaces $\H^n$ which admits cocompact
reflection groups. Thus 2-dimensional examples come from regular
$p$-gons, $p\ge 5$ in the hyperbolic plane. In dimension 3 there
exists a right-angled regular hyperbolic dodecahedron ([A]). In
dimension 4 there exists a right-angled hyperbolic 120-cell
([C],[D2],[PV]). The examples exist up to dimension 8 [VS]. Of
course the dimension of the hyperbolic spaces is limited by
Vinberg's theorem ($\le 29$) [V]. Existence of aperiodic sets of
tiles for these cases also follows from results of Block and
Weinberger [BW]. A new part of our results is the limit strong
aperiodicity of tilings. Also the Block-Weinberger aperiodic
tilings are unbalanced. In this paper we construct strongly
balanced tilings which are limit strongly aperiodic. A strongly
aperiodic tiling of $\H^2$ was recently constructed by
Goodman-Strauss [GS] (his tiling has even finite strongly
aperiodic set of tiles).

First we obtain our tilings of the Davis complex as a tiling by color
with geometrically the same tile (the chamber). Then we take a geometric
resolution of the tiling by color.
If a discrete group $G$ acts by isometries properly cocompactly
on a metric space $X$, there
is an universal way to construct an aperiodic tiling of $X$ by means of an
aperiodic coloring of $G$. We consider an orbit $Gx$ and consider the
Voronoi cells $V_{gx}$ for $g\in G$ where
$$ V_{y}=\{x\in X\mid d(y,x)\le d(y,Gx)\}.$$
Clearly, all cells are isometric to each other and an aperiodic coloring
of $G$ defines an aperiodic tiling of $X$ by color.
In the case of Coxeter groups one can consider colorings of the walls instead of
groups. This allows us to construct strictly
balanced aperiodic tiling.

This paper is arranged as follows. First we consider coloring of discrete
groups (\S1). Then we extend this to a coloring of spaces, in particular trees,
on which the group acts (\S2). Then we apply this to the case of trees of walls in
a Coxeter group (\S3). Finally we construct a strictly balanced limit aperiodic
tilings (\S5). In (\S4) we give a brief account of topology on
the space of tilings.

It is a pleasure to thank Victor Bangert and Mark Sapir for useful discussion
about the Morse-Thue sequence. Also we would like to thank the Max-Plank Institut
f\"ur Mathematik in Bonn for its hospitality.

\head \S1 Aperiodic coloring of groups \endhead

DEFINITION. {\it A coloring of a set $X$ by the set of colors $F$
is a map $\phi:X\to F$.}

We consider the product topology on the set of all colorings $F^X$
of $X$ where $F$ is taken with the discrete topology.

A group $G$ acts (from the left) on the space of functions $F^G$
by the formula $(g\phi)(x)=\phi(g^{-1}x)$.

DEFINITION. A coloring $\phi:\Gamma\to F$ of a discrete group is
called {\it aperiodic} if $\phi\ne g\phi$ for all
$g\in\Gamma\setminus {e}$. It means that the orbit $\Gamma\phi$ of
$\phi$ under the left action of $\Gamma$ on $F^{\Gamma}$ is full
i.e. every the isotropy group of $\phi$ is trivial.

If $\phi=g\phi$ for some $g\in\Gamma$ we call $\phi$ {\it
$g$-periodic}. Clearly, the $g$-periodicity is equivalent to the
$g^{-1}$-periodicity.

Note that every group admits an aperiodic coloring by two elements
$\delta_e:\Gamma\to\{0,1\}$ with $\delta_e(e)=1$ and
$\delta_e(g)=0$ for $g\ne e$. This coloring is not interesting
since it fails to satisfy the following condition.

DEFINITION. A coloring $\phi:\Gamma\to F$ of a discrete group is
called {\it limit aperiodic} if the action of $\Gamma$ on the
closure $\overline{\Gamma\phi}$ of the orbit $\Gamma\phi$ is free,
i.e., every coloring $\psi\in\overline{\Gamma\phi}$ is aperiodic.

We consider the question whether every finitely generated group
admits a limit aperiodic coloring by finitely many colors.

REMARK. There is a weaker version of this question. We say that a
coloring $\phi\in F^{\Gamma}$ is {\it weakly aperiodic} if the
orbit $\Gamma\phi$ is infinite.  A coloring $\phi$  is called {\it
limit weakly aperiodic} if every coloring in the closure
$\psi\in\overline{\Gamma\phi}$ is weakly aperiodic. We note that
the weakly aperiodic version of this question has an affirmative
answer. Namely V. Uspenskii proved [Us] that for every discrete
group $\Gamma$ the topological dynamical system
$(F^{\Gamma},\Gamma)$ has a compact infinite $\Gamma$-invariant
set $X\subset F^{\Gamma}$ with a minimal action on it. We recall
that an action is {\it minimal} if every orbit is dense. Thus, no
orbit in $X$ can be finite and hence every element $\phi\in X$ is
limit weakly aperiodic.

We don't deal with weakly aperiodic colorings in this paper. We
just note that weakly aperiodic colorings correspond to aperiodic
tilings and aperiodic coloring correspond to strongly aperiodic
tilings (see \S 4).

\proclaim{Proposition 1} Let $H\subset G$ be a finite index
subgroup. Then the group $G$ admits a limit aperiodic coloring by
finitely many colors if and only if $H$ does.
\endproclaim
\demo{Proof}
Let $\phi:G\to F$ be a limit aperiodic coloring of $G$. Let $n=[G:H]$
and let $G=\coprod_{i=1}^n Hy_i$ be the decomposition of $G$ into the right
$H$-cosets. Let $\phi y$ denote the result of the right $y$-action, that is
$(\phi y)(x)=\phi(xy)$.

We define $\psi:H\to F^n$ by the formula
$\psi(x)=((\phi y_1)(x),\dots,(\phi y_n)(x))$. Assume that
$\lim_k(h_k\psi)=\psi'$ for a sequence $h_k\in H$ and $\psi'(ax)=\psi'(x)$
for some $a\in H$ and for all
$x\in H$. Taking a subsequence we may assume that the sequence
$h_k\phi$ is convergent. Since
$\lim_k(h_k\phi)$ is not $a$-periodic,
there is $z\in G$ such that $(h_k\phi)(az)\ne(h_k\phi)(z)$ for infinitely many
$k$. Let $z\in Hy_i$. Then for all sufficiently large $k$ we
have the equality $(h_k\psi)(ax)=(h_k\psi)(x)$ for
$x=zy^{-1}_i\in H$. Hence $(h_k\phi y_i)(ax)=(h_k\phi y_i)(x)$. Therefore
we have a contradiction:
$$(h_k\phi)(az)=\phi(h_k^{-1}azy^{-1}_iy_i)=
(h_k\phi y_i)(ax)=(h_k\phi y_i)(x)=
\phi(h_k^{-1}xy_i)=(h_k\phi)(z).$$

In the other direction we may assume that $H$ is normal. If it is not,
we take a smaller
normal subgroup of finite index $H'$. By the above $H'$ admits a
limit aperiodic coloring.
Let $\psi:H\to F$ be a limit aperiodic coloring of $H$ and let $n=[G:H]$.
We define a coloring $\phi:G\to F\times\{1,\dots,n\}$ by the formula
$\phi(x)=(\psi(y_i^{-1}x),i)$ for $x\in Hy_i$ where
$G=\coprod_{i=1}^n Hy_i$ is
the decomposition of $G$ into the left
$H$-cosets (= the right $H$-cosets).
Assume that $\phi'=\lim_kg_k\phi$ is $a$-periodic for some $a$.
Since all $H$-cosets are colored by different colors, $a$ must be in $H$.
We may assume that all $g_k\in y_jH$ for some fixed $j$.
Let $a'=y_j^{-1}ay_j$.
Since $\psi$ is limit aperiodic, there is $z\in H$ such that
$(y^{-1}_jg_k\psi)(z)\ne(y^{-1}_jg_k\psi)(a'z)$. We take $x=y_jz$.
Then for infinitely many $k$, $g_k\phi(x)=g_k\phi(ax)$.
Note that $g_k\phi(x)=\phi(g_k^{-1}y_jz)=\psi(g^{-1}_ky_jz)$ by
the definition of $\psi$ and the choice $y_1=e$.
Thus, $g_k\phi(x)=(y_j^{-1}g_k\psi)(z)$. On the other hand
$g_k\phi(x)=g_k\phi(ax)=\phi(g_k^{-1}ay_jz)=\phi(g_k^{-1}y_ja'z)=
(y_j^{-1}g_k\psi)(a'z)$. Contradiction.
\qed
\enddemo

EXAMPLE. The following coloring of $\Z$ is not limit
aperiodic:
$$
\phi(n)=\cases 1\ \text{if}\ n=\pm k^2, k\in\N\\
               0 \ \text{otherwise}.\endcases
$$
We consider the sequence of colorings $\psi_m(x)=\phi(x+m^2)$,
$\psi_m\in\Z\phi$, and show that  has the
constant $0$-coloring is the limit of $\psi_m$. We need to show that
$\lim_{m\to\infty}\psi_m(x)=0$ for every $x\in\N$. Since the
equation $x+m^2=\pm k^2$ has only finitely many integral solutions
$(m,k)$, the result follows.

Since every weakly aperiodic coloring of $\Z$ is aperiodic, the
existence of a limit aperiodic coloring of $\Z$ follows from the
existence of non-periodic minimal set for the shift action of $\Z$
on $\{0,1\}^{\Z}$ (see the above REMARK).

An explicit example of a limit aperiodic coloring of $\Z$ can be
given by means of the Morse-Thue sequence $m:\N\to\{0,1\}$.

DEFINITION [Mor],[Th]. Morse-Thue sequence $m(i)$. Consider a substitution
rule: $0\to 01$ and $1\to 10$. Then start from $0$ to perform this
substitutions
$$
0\to 01\to 0110\to 01101001\to\dots.
$$
By taking the limit we obtain a sequence of 0 and 1 called the
{\it Morse-Thue sequence}.

\proclaim{Theorem 1 [HM],[Th]} The Morse-Thue sequence contains no
string of type $WWW$ where $W$ is any word in 0 and 1.
\endproclaim
We consider the coloring $\phi:\Z\to\{0,1\}$ defined as
$\phi(x)=m(|x|)$.
\proclaim{Proposition 2} The coloring $\phi$ of $\Z$ is limit
aperiodic.
\endproclaim
\demo{Proof}
Assume the contrary: there is a sequence $\{n_k\}$
tending to infinity and $a\in\N$
such that $\psi(x+a)=\psi(x)$ where $\psi(x)=\lim_{k\to\infty}\phi(x+n_k)$ for
all $x\in\Z$. We may assume that all $n_k>0$.
Then there is $k_0$ such that for all $x\in[1,3a]$ and all $k>k_0$, we have
$\psi(x)=\phi(x+n_k)$. Let $w=\psi(1),\psi(2)\dots,\psi(a)$. Note that
$\psi(1),\psi(2),\dots,\psi(3a)=www$. On the other hand
$\psi(x)=\phi(x+n_k)=m(x+n_k)$ for $x\in[1,3a]$. Thus we have
a "cube" $www$ in the Morse-Thue sequence. Contradiction
\qed
\enddemo

QUESTION: Actually this coloring of $\Z$ has the following
stronger property:
given $n \in \Z \setminus 0$ and $m \in \Z$ there exists a
$q \in \Z$ with
$|q-m| \le 3 |n|$ such that
$\phi(q) \ne \phi(q+n)$.

One can ask, if every finitely generated group has a finite coloring with a
similar property. To be precise, let $G$ be a finitely generated group, and
$d$ be the word metric with respect to a finite generating set.

Does there exists a finite coloring
$\phi:G \to F$ and a constant $\lambda > 0$, such that for every element
$g \in G \setminus e$ and every $h \in G$ there exists
$b \in B_{\lambda d_g(h)}(h)$ with $\phi(gb)\ne \phi(b)$ ?
Here $d_g(h)=d(gh,h)$ is the displacement of $g$ at $h$ and
$B_r(h)$ the distance ball of radius $r$ with center $h$.

Such a coloring $\phi$ is not $g$ invariant and one can see this aperiodicity
already by considering the coloring only on a distance ball $B_r(h)$,
where $h$ is an arbitrary point and the radius
is a fixed constant times the dispacement $d_g(h)$.
A coloring with this property can be considered as a 
natural generalization of
the Morse-Thue coloring to the group G. Such a coloring is 
in some sense
"as aperiodic as possible" and in particular limit aperiodic.

\head \S2 Aperiodic coloring of hyperbolic groups \endhead

It is very plausible that every finitely generated group has limit
aperiodic colorings by finitely many colors. On the other hand a
random coloring is not limit aperiodic. In this section we
construct such colorings for torsion free finitely generated
hyperbolic groups.

In the sequel $G$ is a finitely generated group, and $\|.\|$ the
norm with respect to a finite set of generators.

\proclaim{Lemma 1}
Let $a\in G$ be an element of infinite order in a hyperbolic group.
Then there is $n=n(a)$ such that for every $g\in G$ there is $k\le n$ with
$\|ga^k\|\ne\|g\|$.
\endproclaim
\demo{Proof} It is known that the sequences $\{a^k\}$ and
$\{a^{-k}\}$ define two different points on the Gromov boundary
$\partial_{\infty}G$ of $G$. Let $\xi:\R\to K$ be a geodesic in
the Cayley graph $K$ joining these two points. Note that the
action of $a$ on $K$ leaves these points at infinity invariant. Let
$d=d(e,im\xi)$, then $d(a^k,im\xi)=d(e,a^{-k}im\xi)\le d+\delta$
for every $k$ where $G$ is $\delta$-hyperbolic. The last
inequality follows from the fact that a degenerated triangle in
$K$ defined by the geodesics $im\xi$ and $a^{-k}(im\xi)$ is
$\delta$-thin. Take $n$ such that $\|a^n\|>2\|a\|+10d+14\delta$.
Assume that there is $g$ such that
$\|g\|=\|ga\|=\|ga^2\|=\dots=\|ga^n\|$. Consider the geodesic
$g(im(\xi))$. Let $w\in g(im(\xi))$ and $w'\in g(im(\xi))$ be the
closest points in $g(im(\xi))$ to $g$ and $ga^n$ respectively.

Since the triangle $<e,w,w'>$ is
$\delta$-thin, the geodesic segment $[w,w']$ contains a point $z$ such that
$d(z,y)<\delta$ and $d(z,y')<\delta$ where
$y\in[e,w]$ and $y'\in[e,w']$. Thus,
$$   (1) \ \ \ \ \ \ \ \ \ \ \ \ \ \ \ \ \ \ \ \ \ \ \ \ \ \ \ \ \ \ \ 2\|z\|-\|y\|-\|y'\|\le 2\delta.$$

Denote by $z_k$, $k=0,\dots,n$, a point in $g(im(\xi))$ such that
$d(z_k,ga^k)\le d+\delta$. There is $i$ such that $z\in[z_i,z_{i+1}]\subset
g(im(\xi))$.
Then $d(z,ga^i)\le 3d+3\delta+\|a\|$. Thus,
$$ (2)\ \ \ \ \ \ \ \ \ \ \ \ \ \ \ \ \ \ \ \|g\|=\|ga^i\|\le\|z\|+3d+3\delta+\|a\|.$$

In view of (1) and (2) and the facts $|\|w\|-\|g\||,|\|w'\|-\|g\||\le d+\delta$
we obtain a contradiction:
$\|a^n\|=d(g,ga^n)\le d(w,w')+2d+2\delta=d(w,z)+d(z,w')+2d+2\delta\le
d(w,y)+d(y',w')+2d+4\delta=\|w\|-\|y\|+\|w'\|-\|y'\|+2d+4\delta\le
2\|g\|-\|y\|-\|y'\|+4d+6\delta\le 2\|z\|-\|y\|-\|y'\|+ 2\|a\|+10d+12\delta
\le 2\|a\|+10d+14\delta$.
\qed
\enddemo

EXAMPLE. The group $\Z^2$ does not have the above property with respect to
the generators $(0,\pm 1)$ and $(\pm 1,0)$. Take $a=(1,-1)$ and $g_n=(n,n)$.
Then $\|(n,n)-k(1,-1)\|=2n$ for all $k\le n$.

A geodesic segment in a finitely generated group is the corresponding
sequence of vertices in a geodesic segment in
the Cayley graph. A geodesic segment $[x_1,\dots,x_k]$ is
called {\it radial} if $\|x_1\|<\|x_2\|<\dots<\|x_k\|$.

To construct the limit aperiodic coloring we consider a square
free sequence $\nu:\N\to\{0,1,2\}$, i.e. a sequence which does not
contain any subsequence of the form $WW$, where $W$ is a word in
$0,1,2$. It is possible to construct a square free sequence in the
following way: take the Thue-Morse sequence $0110100110010110...$
and look at the sequence of words of length 2 that appear: $01$,
$11$, $10$, $01$, $10$, $00$, $01$, $11$, $10$ .... Replace $01$
by $0$, $10$ by $1$, $00$ by $2$ and $11$ by $2$  to get the
following: $021012021 \ldots$. Then this sequence is square-free
 [HM].
\proclaim{Theorem 2}
Every finitely generated torsion free hyperbolic
group admits
a limit aperiodic coloring by $9$ colors.
\endproclaim
\demo{Proof} The set of colors will be the set of pairs
$(m,n)$ where $m,n\in\{0,1,2\}$.
Let $G$ be a group taken with the word metric with respect a finite
generating set $S$. We define $\phi(g)=(\nu(\|g\|),\|g\|\ mod\ 3)$
for every $g\in\G$.

Claim : {\it Let $[x_1,x_2,\dots,x_k]$ be
a radial geodesic segment.
Let $g\in G$ be such that
$\phi(x_i)=\phi(g x_i)$
for all $i\in \{1,\ldots k\}$.
Assume in addition that $d(g x_{i_0}, x_{i_0}) < \frac{k}{2}$
for some $i_0$, $1\le i_0\le k$.
Then
$\|g x_i\|=\|x_i\|$ for all $i$.}

We first show that
$[g x_1,\dots,g x_k]$ is also a radial geodesic segment.
Since multiplication from the left is an isometry,
$[g x_1,\dots,g x_k]$ is
a geodesic segment and in particular
$ -1 \leq \|g x_{i+1}\| - \|g x_i\| \leq 1$.
Since $g$ preserves the coloring we have
$\|g x_{i+1}\| \equiv \|g x_i\|+1 \ mod\ 3$.
This two relations imply
$\|g x_{i+1}\|= \|g x_i\| +1$, hence
$[g x_1,\dots,g x_k]$ is a radial geodesic segment
and thus
$\|g x_i\| = \|x_i\| + q$ for some fixed integer $q$.
By our assumption $|q| < \frac{k}{2}$.
Assume $0 < q$. Since $g$ preserves the colors we obtain
the equality
$$(\nu(\|x_1\|),\dots,\nu(\|x_{q}\|))=(\nu(\|x_{q+1}\|),\dots,\nu(\|x_{2q}\|)),$$
a contradiction to the square freeness of $\nu$.
In a similar way we obtain a contradiction if $q < 0$.
This implies $q=0$ and hence the claim.

Now assume that there is a sequence $g_l\in G$ with
$\|g_l\|\to\infty$ such that the limit
$\lim_{l\to\infty}\phi(g_lx)=\psi(x)$ exists for every $x\in G$.

Assume that there is $b\in G\setminus \{e\}$ such that $\psi(x)=
\psi(bx)$ for all $x$.

Let $n$ be taken from Lemma 1 for $a=b$.
We may assume that there is $l_0$ such that
for $l>l_0$, $\psi(y)=\phi(g_ly)$ for all $y\in B_{4n\|b\|}(e)$.
Consider a radial geodesic segment
$x_{1},\dots ,x_{k}$ of the length $k-1$ with $k=3n\|b\|$
and with endpoint $x_k=g_l$. Such a segment clearly exists for all $l$
large enough.

Let $s$ be the smallest natural number such that
$\|g_lb^s\|\ne\|g_l\|$, thus $s \le n$.
Both segments $g_l^{-1}([x_1,\dots x_{k}])$ and
$g_l^{-1}(g_lb^sg_l^{-1}[x_{1},\dots,x_{k}])$
lie in $B_{4n\|b\|}(e)$.
Thus
$\phi(x)=\phi(g_lg_l^{-1}(x))
=\psi(g_l^{-1}(x))=\psi(b^sg_l^{-1}(x))=
\phi(g_lb^sg_l^{-1}x)$ for all $x\in[x_1,\dots ,x_{k}]$.
Furthermore we compute
$d(g_lb^sg_l^{-1}x_k,x_k)=d(g_lb^s,g_l)=\|b^s\|\le n\|b\|<\frac{k}{2}$.
We apply
the claim to $g=g_lb^sg_l^{-1}$ to obtain
the contradiction: $\|g_lb^s\|=
\|g_lb^sg_l^{-1}x_k\|=\|x_k\|=\|g_l\|$.
\qed
\enddemo

REMARK. It is still an open problem whether every hyperbolic group
contains a torsion free subgroup of finite index.

\head \S3 Aperiodic coloring of $G$-spaces \endhead

We note that in the definitions from the beginning of \S1 one can
replace a group $G$ by a space $X$ with a $G$-action. Thus, $G$
acts on the space of colorings $F^X$ by the same formula
$(g\phi)(x)=\phi(g^{-1}x)$. Let $K$ be the kernel of the action.
We can speak about {\it $G$-aperiodic colorings} of $X$ as ones
with $\phi\ne g\phi$ for all $g\in G\setminus K$. Similarly one
can define {\it limit $G$-aperiodic} colorings of $X$ as those
whose orbit $G\phi$ has only $G$-aperiodic coloring in its closure
$\overline{G\phi}$ with respect to the product topology $F^X$.

The following is an analog of Proposition 1.
\proclaim{Proposition 3}
Let $H\subset G$ be a finite index subgroup and let $G$ act on $X$.
Then $X$ admits a limit $G$-aperiodic coloring by finitely many colors if
and only if it admits a limit $H$-aperiodic coloring by finitely many colors.
\endproclaim

It is possible to extend Theorem 2 to the case of an isometric
action on a hyperbolic space $X$.
We consider only the case when $X$ is a simplicial tree. Thus
every edge has length equal to $1$.

If $x_0$ is a root of $X$, we denote by $\|x\|=d(x,x_0)$ for $x\in
X$. We prove the following analog of Lemma 1.

\proclaim{Lemma 2} Let $G$ act on a simplicial rooted tree $X$ and
let $a\in G$ operate without fixed points. Then  for every $g\in
G$  there is $k\le 2$ with
$\|ga^kx_0\|\ne\|gx_0\|$.
\endproclaim
\demo{Proof} Assume that $\|ga^2x_0\|=\|gax_0\|=\|gx_0\|$. Since
$d(gx_0,gax_0)=d(hgx_0,hgax_0)=d(gax_0,ga^2x_0)$ for $h=gag^{-1}$,
the points $z=gx_0$, $h(z)=gax_0$, and $h^2(z)=ga^2x_0$ have common predecessor
$y$ in the tree which is the common midpoint of the geodesic segments
$[z,h(z)]$ and $[h(z),h^2(z)]$. Then $hy=y$ and hence $g^{-1}y$ i
s a fixed point for $a$:
$a(g^{-1}y)=g^{-1}y$.
This contradicts to the assumption.

\qed
\enddemo

Let $x_0\in X$ be a base point in a tree $X$. We consider the
coloring of the set of vertices of $X$ defined like in the proof
of Theorem 2: $\phi(x)=(\nu(\|x\|),\|x\|\ mod\ 3)$.
\proclaim{Proposition 4} Suppose that a group $G$ acts by
isometries on a simplicial tree $X$ with the above coloring $\phi$
and let $\psi\in\overline{G\phi}$. Then $b\psi\ne\psi$ for every
$b\in G$ with unbounded orbit $\{b^kx_0\mid k\in N\}$. Moreover,
$b\psi\ne\psi$ on the orbit $Gx_0$.
\endproclaim
\demo{Proof}
First we note that a similar as in the proof of Theorem 2
claim take place:

{\it Let $[z_1,\dots,z_{k}]$ be a radial geodesic segment in $X$.
Let $g \in G$ such that $\phi(g z_i) = \phi(z_i)$ for all $i$ and
$d(g z_{i_0}, z_{i_0}) < \frac{k}{2}$ for some $i_0$, $1\le i_0\le
k$. Then $\|g z_i\|=\|z_i\|$.}

The proof is exactly the same as in Theorem 2.

Assume that there is a sequence $g_l\in G$ with
$\|g_l\|\to\infty$ with the limit $g_l^{-1}\phi$ equal
to a $b$-periodic coloring $\psi$ such that $\{b^kx_0\}$ is infinite.
That is the limit
$\lim_{k\to\infty}\phi(g_kx)=\psi(x)$ exists for every $x\in X$.

Since the orbit $b^mx_0$ is infinite we can apply Lemma 2  for $a=b$. 
Let $[y_1,\dots,  y_{8k}]$ be a radial
segment with $y_1=x_0$ where $k=6\|b(x_0)\|$. By the definition
of the point-wise limit we may assume that there is $l_0$ such
that for $l>l_0$, $\psi(y)=\phi(g_ly)$ for all $y\in [y_1\dots y_{8k}]
\cup b[y_1\dots y_{8k}]
\cup b^2[y_1\dots y_{8k}]$. Then the image $g_l[y_1,\dots
y_{8k}]$ contains either a radial segment $[z_{1},\dots,z_{k}]$ of
the length $k-1$ with $z_k=g_l(x_0)$  or it contains a radial
segment $[z_{1},\dots,z_{6k}]$ with $d(z_1,g_l(x_0))<k$.

Let $s$ be the smallest natural number such that $\|g_lb^sx_0\|\ne\|g_lx_0\|$.
Thus, $s\le 2$.

Then
$\phi(x)=\phi(g_lg_l^{-1}(x))=\psi(g_l^{-1}(x))=\psi(b^sg_l^{-1}(x))=
\phi(g_lb^sg_l^{-1}x)$ for all $x\in[z_{1},\dots ,z_{k}]$.
Furthermore, in the first case we compute
$d(g_lb^sg_l^{-1}z_k,z_k)=d(b^sx_0,x_0)=\|b^sx_0\|\le 2\|bx_0\|<
\frac{k}{2}$. We apply the claim to $g = g_lb^sg_l^{-1}$  with
$i_0=k$ to obtain the contradiction: $\|g_lb^sx_0\|=
\|g_lb^sg_l^{-1}z_k\|=\|z_k\|=\|g_lx_0\|$.

In the second case $d(z_1,g_lb^sg_l^{-1}z_1)\le$
$$d(z_1,g_l(x_0))+
d(g_l(x_0),g_lb^s(x_0))+d(g_lb^s(x_0),g_lb^sg_l^{-1}z_1)\le
2k+2\|b(x_0)\|\le 3k.$$ We apply the claim with $g =
g_lb^sg_l^{-1}$, $i_0=1$, and $6k$ instead of $k$. Let $i$ be such
that $g_ly_i=z_1$. Then $\|z_i\|=\|g_lx_0\|$. Since
$\phi(\|g_lb^sy_{i-1}\|)=\phi(\|g_ly_{i-1}\|)\ne\phi(\|g_ly_i\|-1)$,
we obtain that $\|g_lb^sy_{i-1}\|\ne\|g_ly_i\|-1=\|z_1\|-1$. Hence
$\|g_lb^sy_{i-1}\|=\|z_1\|+1$ and $\|g_lb^sx_0\|=\|z_1\|+i-1$. The
we obtain a contradiction: $\|g_lb^sx_0\|=\|z_1\|+i-1=
\|g_lb^sg_l^{-1}z_1\|+i-1=\|g_lb^sg_l^{-1}z_i\|=\|z_i\|=\|g_lx_0\|$
. \qed
\enddemo
As a consequence we obtain the following:
\proclaim{Corollary 1}
Suppose that a group $G$ acts on a rooted simplicial tree
$X$ such that $Gx_0$ is a full orbit. Then
there is a limit $G$-aperiodic coloring of $X$ by 9 colors.
Moreover, the restriction of this coloring to the orbit $Gx_0$ is
also limit $G$-aperiodic.
\endproclaim

\proclaim{Theorem 3} Suppose that a group $G$ acts by isometries
on  simplicial trees $X_1,\dots, X_n$ in
such a way that the induced action on the
product $\prod X_i$ is free. Then $G$ admits a limit aperiodic coloring
by $9n$ colors.
\endproclaim
\demo{Proof}
First we note that the fixed point theorem for CAT(0) spaces
implies that $G$ must be torsion free.

Let $K_i$ be the kernel of the representation $G\to Aut(X_i)$ and let
$G_i=G/K_i$.
We fix a base point $x_0^i$ in each tree
and consider  colorings
$\phi^i:X_i\to F_i$, $|F_i|=9$ from Proposition 4. This defines a map
$\phi:\prod X_i\to\prod F_i$. Let $\phi':G\to\prod F_i$ be the restriction
to the orbit: $\phi'(g)=\phi(gx_0)$ where $x_0=(x^1_0,\dots,x^n_0)$.

Assume that $\psi=\lim g^{-1}_k\phi'$ be $a$-periodic:
$a\psi(x)=\psi(x)$, $a\in G\setminus\{e\}$. Then
$\psi=(\psi_1,\dots,\psi_n)$ where $\psi_i=\lim g^{-1}_k\phi'_i$
and $\phi'_i=\phi_i|_{Gx_0^i}$. Since $a$ operates without fixed
points on $X$, there exists $i$ such that $a$ has no fixed points
on $X_i$. Then by Proposition 4 $a\psi\ne\psi$ which contradicts
with the assumption. \qed
\enddemo

We recall the definition of Coxeter groups
A {\it Coxeter matrix} $M=(m_{ij})$ is a
symmetric square matrix with 1 on the diagonal and all other entries
are from $\N_+=\{0\}\cup\N$.
A {\it  Coxeter group} $\Gamma$ with a generating set $S$
is a group with a presentation
$$
\langle S\mid (uv)^{m_{ij}}=1, m_{ij}\in M\rangle
$$
where $M$ is a  Coxeter matrix. Here we use the convention $a^0=1$.
Traditionally the literature on Coxeter groups uses $\infty$ instead of $0$.

A Coxeter group is called {\it right-angled} if all $m_{ij}$ are $0$s or $2$s.

\proclaim{Theorem 4} Every Coxeter group $\Gamma$
admits a limit aperiodic coloring by finitely many colors.
\endproclaim
\demo{Proof}
Since every Coxeter group contains a finite index torsion free
subgroup, in view of Proposition 1 it suffices to prove it for
a finite index subgroup.
We apply Januszkiewicz's construction of equivariant isometric embedding of
a torsion free finite index normal subgroup $\Gamma'\subset\Gamma$
into the finite product of trees [DJ] and then apply Theorem 3.\qed
\enddemo

Let $K$ be the Cayley graph of a Coxeter group $(\Gamma,S)$.
Then every generator $c\in S$ every conjugate $w=gcg^{-1}$, $g\in\Gamma$,
acts on $K$ by reflection.
Let $M_w$ be the set of fixed point of $w$. We call it the {\it
wall (or mirror)} of the reflection $w$. Clearly, that $uM_w$ is a wall
for all $u$ and $w$. Therefore $\Gamma$ acts on the set of all walls $\Cal W$.
According to [DJ] all walls can be partitioned in finitely many classes
$\Cal W=\Cal W_1\cup\dots\cup\Cal W_m$
such that the walls from each class $\Cal W_i$
form a vertex set of a simplicial tree $T_i$. Moreover, all sets $\Cal W_i$
are invariant under
a normal finite index subgroup $\Gamma'$ which acts by isometries on each $T_i$
in such a way that the $\Gamma'$-action on the product $\prod_{i=1}^m T_i$
is free.
\proclaim{Theorem 5}
The set $\Cal W$ of all walls in a Coxeter group admits a limit
$\Gamma$-aperiodic
coloring by finitely many colors.
\endproclaim
\demo{Proof}
On every tree $T_i$ we consider a coloring from Proposition 6 with 9 each time
different colors. Thus we use $9m$ colors. This defines a coloring
$\phi=\cup\phi_i$ of $\Cal W$.
Let $\psi$ be a limit coloring of $\Cal W$. Then, clearly, $\psi=\cup\psi_i$
where each $\psi_i$ is a limit coloring for $\phi_i$. Let $b\in\Gamma'$.
Since the action of $\Gamma'$ on the product $\prod_{i=1}^m T_i$ is free
and $b$ is of infinite order, we obtain that $b$ has an infinite orbit on
some $T_i$. By Proposition 6 $b\psi_i\ne\psi_i$ and hence $b\psi\ne\psi$.
Proposition 5 completes the proof.
\qed
\enddemo

\head \S4  Space of tilings
\endhead

We recall a definition of tiling of a metric space from [BW]. Let
$X$ be a metric space. A set of tiles $(\Cal T,\Cal F)$ is a
finite collection of $n$-dimensional complexes $t\in\Cal T$ and a
collection of subcomplexes $f\in\Cal F$ of dimension $< n$,
together with an opposition function $o:\Cal F\to\Cal F$,
$o^2=id$. A space $X$ is tiled by the set $(\Cal T,\Cal F)$ if
\roster \item{} $X=\cup_{\lambda}t_{\lambda}$ where each
$t_{\lambda}$ is isometric to one of the tiles in $\Cal F$;
\item{} $t_{\lambda}\setminus\cup_{f\in
t_{\lambda}}=Int(t_{\lambda})$ in $X$ for every $\lambda$; \item{}
If $Int(t_{\lambda}\cup t_{\lambda'})\ne Int(t_{\lambda})\cup
Int(t_{\lambda'})$ then $t_{\lambda}$ and $t_{\lambda'}$ intersect
along $f\in t_{\lambda}$ and $o(f)\in t_{\lambda'}$; \item{} There
are no free faces of $t_{\lambda}$.
\endroster

Strictly speaking a tiling of $X$ is a collection $\alpha=\{\phi_{\lambda}\}$
of isometries $\phi_{\lambda}:t_{\lambda}\to t$, $t\in\Cal T$, satisfying the
above axioms. For every tiling $\alpha$ 
there is a minimal (or reduced) set of tiles
$(\Cal T_{\alpha},\Cal F_{\alpha})\subset(\Cal T,\Cal F)$.

Let $X$ be a metric space with a base point $x_0$ Assume that
$\diam t\le 1$ for all $t\in\Cal T$. Let $\alpha$ be a tiling
$X=\cup_{\lambda}t_{\lambda}$ of $X$. We denote by
$\alpha_n=\{t_{\lambda}\mid t_{\lambda}\subset B_n(x_0)\}$ where
$B_r(x)$ stands for the closed ball of radius $r$ centered at $x$.
For a metric space $Y$ we denote by $exp Y$ the space of compact
subsets of $Y$ taken with the Hausdorff metric. Also for $m\in\N$
denote by $exp_mY$ the $m$th hyperpower of $Y$, i.e., subspace of
$exp Y$ that consists of subsets of cardinality $\le m$. Note that
$\alpha_n$ defines a point in $exp(exp B_n(x_0))$. Actually, there
is $k=k(n)$ such that $\alpha_n$ lies in $\exp_k(exp B_n(x_0))$.
Clearly, the sequence $\alpha_n$ completely defines the tiling
$\alpha$.

The space of tilings was defined by many authors (see for example
[BBG],[S],[SW]). Here we give an alternative definition.
Let $tl(X,\Cal T)$ denote the set of all $\Cal T$-tilings of $X$.
We introduce the topology on $tl(X,\Cal T)$
as a subspace topology
$$tl(X,\Cal T)\subset\prod_{n=1}^{\infty} exp(exp B_n(x_0)).$$
Let $\Cal F$ be a finite family of compact subsets in a metric space $Y$.
We denote by $exp^{\Cal F}Y$ the subspace of $exp Y$ whose points are isometric
copies of elements of $\Cal F$.

Note that $tl(X,\Cal T)\subset\prod_{n=1}^{\infty}exp_{k(n)}exp^{\Cal T}(B_n)$.

The following Proposition is well-known [S].
\proclaim{Proposition 7} The space $tl(X,\Cal T)$ is compact.
\endproclaim
\demo{Proof}
Since $exp^{\Cal T}(B_n)$ is compact, it suffices to show that $tl(X,\Cal T)$
is a closed subset in $\prod_{n=1}^{\infty}exp_{k(n)}exp^{\Cal T}(B_n)$.
For that it suffices to show that the set
$\{\alpha_n\mid\alpha\in tl(X,\Cal T)\}$ is closed in
$exp_{k(n)}exp^{\Cal T}(B_n)$ for every $n$.
\qed
\enddemo

Let $G\subset Iso(X)$ be a subgroup of the group of isometries of $X$.
Clearly $G$ acts on
$tl(X,\Cal T)$. We say that a tiling $\alpha\in tl(X,\Cal T)$ is
{\it strongly $G$-aperiodic} if $g\alpha\ne\alpha$ for all $g\in G\setminus\{e\}$.
A tiling $\alpha$ is called {\it aperiodic} if the group $Iso(\alpha)\subset
Iso(X)$ of isometries
of $\alpha$ does not
act cocompactly on $X$. A tiling $\alpha$ is {\it limit strongly $G$-aperiodic} if
every tiling $\beta\in \overline{G\alpha}$ is strongly aperiodic.
If $G=Iso(X)$ we use the terms {\it strongly aperiodic} and {\it limit strongly
aperiodic}.

\head \S5  Aperiodic tiling of Davis complex \endhead

Here we recall the definition of the Davis complex [D1]. Let
$\Gamma$ be a Coxeter group with generating set $S$. The nerve
$N=N(\Gamma,S)$ is the simplicial complex defined in the following
way: the vertices of $N$ are elements of $S$. Different vertices
$s_1,\dots s_k$ span a simplex $\sigma$ if and only if the set
$s_1,\dots s_k$ generates a finite subgroup $\Gamma_{\sigma}$ of
$\Gamma$. By $N'$ we denote the barycentric subdivision of $N$.
The cone $C=Cone N'$ over $N'$ is called a {\it chamber} for
$\Gamma$. The Davis complex  $X=X(\Gamma,S)$ is the image of a
simplicial map $q:\Gamma\times C\to X$ defined by the following
equivalence relation on the vertices: $a\times v_{\sigma}\sim
b\times v_{\sigma}$ provided $a^{-1}b\in\Gamma_{\sigma}$ where
$\sigma$ is a simplex in $N$ and $v_{\sigma}$ is the barycenter of
$\sigma$. We identify $C$ with the image $q(e\times C)$. The group
$\Gamma$ acts simplicially on $X$ with the orbit space equivalent
to the chamber. Thus, the Davis complex is obtained by gluing the
chambers $\gamma C$, $\gamma\in\Gamma$ along their boundaries.
Note that $X$ admits an equivariant cell structure with the
vertices $X^{(0)}$ equal the cone points of the chambers and with
the 1-skeleton $X^{(1)}$ isomorphic to the Cayley graph of
$\Gamma$. A conjugate $r=wsw^{-1}$ of every generator $s\in S$ is
a reflection. The fixed point set $M_r$ of a reflection $r$ is
called the wall of $r$. Note that walls defined in \S 3 are
obtained from the walls in Davis' complex by the restriction to
the Cayley graph.

\proclaim{Proposition 8}
Every finite coloring $\phi:\Cal W\to F$ of the set of walls of the Davis
complex $X$ defines a tiling $\bar\phi$ of $X$ with $o(f)=f$.
\endproclaim
\demo{Proof}
The set of tiles $\Cal T$ of $\bar\phi$ is the set of chambers with all possible
colorings of their faces. The set of faces $\Cal F$ is the set of all
possible colored faces of the chambers. Set $o(f)=f$. Then all conditions hold.
\qed
\enddemo
We call the tiling $\bar\phi$ as a {\it tiling by coloring} $\phi$.

Let $(\Cal T,\Cal F)$ be a set of tiles. A function $w:\Cal F\to\Z$ is
called {\it a weight function} if $w(o(f))=-w(f)$ for every $f\in\Cal F$.
We recall a definition from [BW].
\proclaim{Definition} A finite set of tiles $(\Cal T,\Cal F)$ is
unbalanced if there is a weight function $w$ such that $\sum_{f\in t}w(f)>0$
for all $t\in\Cal T$.

It is called semibalanced if $\sum_{f\in t}w(f)\ge 0$ for all $t\in\Cal T$.
\endproclaim
We call a set of tiles {\it strictly balanced} if for every nontrivial weight
function $w$ there are tiles $t_+$ and $t_-$ such that
$\sum_{f\in t_+}w(f)>0$ and $\sum_{f\in t_-}w(f)<0$.

A tiling is called {\it strictly balanced} ({\it unbalanced})
if its minimal set of tiles is strictly balanced (unbalanced).

We now associate to every wall in the Davis complex an orientation.
A wall divides the Davis complex into two components.
Roughly speaking the orientation says which of the components is left and which
is right. Let $\bar\phi$ be tiling of the Davis complex $X$ by coloring of
the walls $\phi$ with the set of tiles $(\Cal T,\Cal F)$.
The orientation of the walls define a new tiling $\phi'$ of $X$ with
the set of tiles $(\Cal T',\Cal F')$ where $\Cal F'=\Cal F_+\cup\Cal F_-$
where $\Cal F_+$ and $\Cal F_-$ are copies of $\Cal F$. The face $f\in
t_{\lambda}$ has sign +, if $Int(t_{\lambda})$ is left of the wall and sign -,
if $Int(t_{\lambda})$ is right of the wall. The opposition function
$o:\Cal F'\to\Cal F'$ maps $f_+$ to $f_-$. We call such tiling a
{\it geometric resolution} of a tiling by coloring. This new tiling is not any
more a tiling by coloring. A geometric meaning of this resolution is
that we deform all faces of
a given color and a given sign in the same direction by the same pattern.
For the faces of the same color but of opposite sign we take an opposite
deformation.

The following is obvious.
\proclaim{Lemma 3} Assume that a coloring $\phi:\Cal W\to F$ is limit aperiodic.
Then the tiling by color $\bar\phi$ as well as it's any geometric resolution
is limit strongly aperiodic.
\endproclaim

Note that in the Davis complex every wall has a canonical orientation, by
deciding that the base chamber $C$ is in the left component. Thus we can
indicate the chosen orientation itself by a sign. A wall gets the sign $+$, if
the orientation of the wall is the canonical one and $-$ otherwise.

In [BW] unbalanced tilings of some nonamenable spaces are
constructed. In particular all hyperbolic Coxeter groups admit
such tilings. We can derive this fact using geometric resolutions.
\proclaim{Proposition 9} Every coloring of the walls for a
hyperbolic Coxeter group admits an unbalanced geometric
resolution.
\endproclaim
\demo{Proof} We assign $+$ to every wall. The hyperbolicity
implies that for every chamber $C'$ the numbers of faces of $C'$
whose walls separate $C'$ from the base chamber $C$ is strictly
less than the number of faces whose walls do not separate $C'$ and
$C$. Then for every chamber $C'$ the faces whose walls do not
separate $C'$ from $C$ obtain sign +, all other -. We define a
weight function by sending a positive face to $+1$ and a negative
face to $-1$. \qed
\enddemo
We note that every unbalanced tiling is aperiodic. This fact can be derived
formally from proposition 4.1 [BW]. Since the proof there has some omissions we
present a proof below.
\proclaim{Proposition 10}
Let $(\Cal T,\Cal F)$ be the set of tiles of a geometric realization of a tiling
by coloring of the Davis complex $X$ of a Coxeter group $\Gamma$. Suppose that
the set of tiles $(\Cal T,\Cal F)$ is unbalanced. Then any $(\Cal T,\Cal
F)$-tiling $\alpha$ is aperiodic.
\endproclaim
\demo{Proof}
Let $G$ be a group of isometries of $\alpha$. Then $G\subset\Gamma$. Hence
$G$ is a matrix group. By Selberg Lemma it contains a torsion free subgroup
$G'$ of finite index. Then the orbit space $X/G'$ is compact and admits a
$(\Cal T,\Cal F)$-tiling (Note that by taking $X/G$ as in [BW] we cannot
always obtain a tiling because of free faces). Then we obtain a contradiction:
$$
0<\sum_{t\in X/G'}\sum_{f\in t}w(f)=\sum_{f\in X/G'}(w(f)+w(o(f)))=0.
$$
\qed
\enddemo

\proclaim{Theorem 6}
For every Coxeter group $\Gamma$ for every coloring $\phi:\Cal W\to F$
with the property that walls of the same
color do not intersect, there is
a strictly balanced geometric resolution. Additionally, every
limit tiling of this resolution is strictly balanced.
\endproclaim
\demo{Proof} First we construct a strictly balanced geometric
resolution of $\phi$. Consider the set of walls $\Cal
W_c=\phi^{-1}(c)$ of the same color $c\in F$. Since the walls of
the same color do not intersect, they are ordered by level from
the base chamber. (The level $lev$ is defined by induction. If one
removes the walls $\Cal W_c$ from $X$, the space is divided into
components. Walls from $\Cal W_c$ that bound the component of the
base chamber are of level one. Then drop the walls of level one
and repeat the procedure to get new walls of level one and call
them of level two and so on). We give the walls $\Cal W_c$ signs
in an alternate fashion by the level $(-1)^{lev(M)}$:\
$-+-+-+-+-\dots $.

We show that this geometric resolution is strictly balanced. Let
$w:F_+\cup F_-\to\Z$ be a nontrivial weight function with $w(f_+)=-w(f_-)$.
We show that there are chambers $C_+$ and $C_-$ such that
$$
\sum_{f\in C_+}w(f)>0\ \ \ \ \text{and}\ \ \ \ \ \sum_{f\in C_-}w(f)<0.
$$
Because of the symmetry it suffices to show the first. Since $w$ is nontrivial,
there exists a face $f^0$ which is the common face of two adjacent tiles
$t_{\lambda}$ and $t_{\lambda'}$ such that $w(f^0)\ne 0$. Let $M_0$ denote the
wall that contains $f^0$. Now there are four cases corresponding to the parity
of
the sign of $w(f^0_+)$ and the sign of $M_0$. We discuss only one, and to be
fair not the easiest of the cases: $w(f^0_+)>0$ and the orientation of $M_0$ is
negative.

We assume that $f^0$ in $t_{\lambda}$ has sign - and in $t_{\lambda'}$ has
sign +. We take a number $k$ larger that the number of walls separating
$t_{\lambda'}$ and the base chamber $C$. Let $c\in F$ be a color.
We call $c$ {\it even} if $w(c_+)>0$, {\it odd} if $w(c_+)<0$, and
{\it neutral} if $w(c_+)=0$. We define

$\Cal W^{2k}_{ev}=\{M\in \Cal W\mid \phi(M) \text{is even and}\ lev(M)=2k\}$

$\Cal W^{2k+1}_{odd}=\{M\in \Cal W\mid \phi(M) \text{is even and}\ lev(M)=2k+1\}$

$\Cal W^{k+1}_{0}=\{M\in \Cal W\mid \phi(M) \text{is neutral and}\
lev(M)=k+1\}$.

Claim 1: {\it The set of walls $\Cal W^{2k}_{ev}\cup \Cal W^{2k+1}_{odd}
\cup\Cal W^{k+1}_{0}\cup \{M_0\}$ bounds a bounded set $D$ containing the
chamber $t_{\lambda'}$.}

Clearly, it bounds a convex set in the
Hadamard space $X$. If it is unbounded, then there is a geodesic
ray from $t_{\lambda'}$
to the visual boundary which does not intersect any of our mirrors.
Since we have only finitely many colors, there is a color $c$ such that this ray
intersects infinitely many walls of this color. By the choice of $k$ the first
of this crossed walls has level $\le k+1$. To get to infinity the ray must
cross walls of color $c$ with all levels $\ge k=1$. Thus, one of the
intersected walls is contained in our set.

Claim 2: {If $f$ occurs as a face of a tile $t_{\mu}\subset D$
such that $f\subset\partial D$, then $w(f)>0$.}

We consider cases.

(i) If $M_f=M_0$, then $t_{\mu}$ lies on the same side of $M_0$ as
$t_{\lambda'}$. Then $w(f)=w(f^0_+)>0$.

(ii) If $f$ is a neutral face, then $w(f)=0$ anyway.

(iii) If $f$ is of even color. Then $f$ is contained in a wall $M$ from
$\Cal
W^{2k}_{ev}$ then $t_{\mu}$ lies on the same
side of the wall $M_f$ as the base chamber. Since $M$ has orientation $+$.
Hence $f$ as a face of $t_{\mu}$ gets the sign +. Hence $w(f)=w(f_+)>0$.

(iv) A similar argument applies for $f$, if $f$ is of an odd color.

According to the Claim 1 we have $D=\cup_{i=1}^kC_i$ where $C_1,\dots,C_k$ is
a finite collection of chambers. Then
$$
\sum_{i=1}^k\sum_{f'\in C_i}w(f')=\sum_{f'\in\partial D}w(f')\ge 0
$$
by Claim 2. Since $f^0_+$ is in the last set of faces, we see that
the expression is indeed $>0$. Therefore, $\sum_{f'\in
C_i}w(f')>0$ for some $i$.

This finishes the proof of the first step. Thus we have constructed a
strictly balanced geometric resolution of $\phi$.

Actually the proof of the first step shows more: If we chose for any given color
$c$ an orientations of walls $\Cal W_c$ in the alternate way $+-+-+-...$
or $-+-+-+...$ (and maybe for different colors in a different way), then the
resulting geometric resolution is strictly balanced. Let us call such choice of
orientations as {\it allowed}. The levels of walls depend on the base chamber.
If we define levels with respect to a different chamber, all parities of the
levels will be either preserved or changed to opposite. as a consequence we
obtain the following: if the orientation of tiling by coloring $\phi$ is
allowed, then also the orientation of the tiling by coloring $g\phi$ is allowed
for every $g\in\Gamma$. Thus also all limit tilings of tiling constructed in the
step 1 are strictly balanced.
\qed
\enddemo
\proclaim{Corollary 2}
For every Coxeter group $\Gamma$ there is a strictly balanced strictly aperiodic tiling of
the Davis complex such that every limit tiling is strictly balanced
and strictly aperiodic.
\endproclaim
\demo{Proof}
We apply Theorem 6 to a coloring from Theorem 5. \qed
\enddemo
Note that in the proof of Corollary 2 we used that $\Gamma$ is
the isometry group of
the Davis complex.

In 2-dimensional jigsaw tiling puzzles a geometric resolution is
usually realized by adding rounded tabs out on the sides of the
pieces with corresponding blank cut into intervening sides to
receive the tabs of adjacent pieces. This procedure destroys the
convexity of the pieces. We show that in the case of the
hyperbolic plane $\H^2$ we can modify this construction to obtain
aperiodic and strictly balanced tiling with convex tiles. Compare
also the papers [MM],[Moz]. \proclaim{Theorem 7} (1) For every
$n\ge 3$ there is a strictly balanced limit strongly aperiodic
tiling of $\H^2$ by convex $2n$-gons.

(2) For every $n\ge 3$ there is a finite set of tiles $(\Cal T,\Cal F)$ that
consists of
convex $2n$-gons with limit strongly aperiodic tiling of $\H^2$ such that every
$(\Cal T,\Cal F)$-tiling of $\H^2$ is aperiodic.
\endproclaim
\demo{Proof} (1) Identify $\H^2$ with the Davis complex for the
right-angled Coxeter group $\Gamma$ generated by reflections at a
regular right-angled $2n$-gon. Coloring the sides of the $2n$-gon
in two colors $a$ and $b$ in an alternating fasion induces a
coloring of the walls $\psi:\Cal W\to\{a,b\}$ such that the walls
of the same color do not intersect. The walls of the same color
$c$ define a tree $T_c$ with an action of $\Gamma$ on it such that
the induced $\Gamma$-action on the product $T_a\times T_b$ is free
(see [BDS] or [DJ]). By Theorem 3 we can refine $\psi$ to a limit
aperiodic coloring $\psi:\Cal W\to\{a_i,b_i\}_{i=1,\dots,9}$ We
apply Theorem 6 to obtain limit strictly balanced geometric
resolution $\phi'$. By Lemma 3 it will be limit strongly
$\Gamma$-aperiodic tiling. It's easy to see that the tiling is
limit strongly aperiodic with respect to entire isometry group of
$\H^2$.

Now we define a modification of the tiling $\phi'$. Consider a
vertex of a translate of the $2n$-gon. It is the intersection
point of an $a_i$-wall with a $b_j$-wall. Denote it by $O_{ij}$.
The orientations on these walls define local coordinate system. We
move the vertex $O_{ij}$ by small amount using these coordinates.
We chose a small different numbers $d_{ij}$, $i,j\in\{1,\dots,9\}$
and move $O_{ij}$ to the distance $d_{ij}$ in the direction of the
diagonal of the positive quadrant. After this deformation we
obtain a finite number of new convex tiles, which (for generic
deformations) only allow tilings of $\H^2$ compatible with the
matching rule defined by $\phi'$. New tiling has all desired
properties.

(2) We take the above coloring $\phi$ of the walls and take a geometric
resolution from Proposition 9. Apply Proposition 10 to complete the proof.
\qed
\enddemo
An interesting question is under what conditions the set of tiles
$(\Cal T,\Cal F)$ of a geometric resolution $\phi'$ of a tiling by coloring
of the Davis complex $X$
is (strongly) aperiodic. Clearly, it is strongly aperiodic whenever its orbit
$\Gamma\phi'$ in $tl(X,\Cal T)$ is dense.

\enddocument